\documentclass[12pt,twosided,reqno]{amsart}
\usepackage{amsmath}
\usepackage{amsthm}
\usepackage{amsfonts}
\usepackage{amssymb}

\theoremstyle{theorem}

\theoremstyle{definition}

\def\<{\left < }
\def\>{\right >}
\def\({\left ( }
\def\){\right )}

\def\k{\kappa}

\begin{document}

\vskip .5cm
\title{Centroid of triangles associated with a curve}
 \vskip0.5cm
\thanks{
         \newline\indent { The first author was   supported by Basic Science Research Program through
    the National Research Foundation of Korea (NRF) funded by the Ministry of Education, Science and Technology (2010-0022926).
    }
}

\vskip 0.5cm

\maketitle

\vskip 0.5cm
\centerline{\scshape
  Dong-Soo Kim  and Dong Seo Kim } \vskip .2in

\begin{abstract}
Archimedes showed that the area between a parabola and any chord $AB$
on the parabola is four thirds of the area of triangle $\Delta ABP$, where P is the
point on the parabola at which the tangent is parallel to the chord $AB$. Recently,
this property of parabolas was  proved to be a characteristic property of parabolas.
With the aid of  this characterization of parabolas, using centroid of
triangles associated with a   curve
we present two conditions  which are necessary and sufficient for
a strictly locally convex  curve in the plane to be a parabola.

\end{abstract}

\vskip 1cm

\date{}
\maketitle

\section{Introduction}
 \vskip 0.50cm

We study   strictly locally convex plane curves.
Recall that a regular plane curve $X:I\rightarrow {\mathbb R}^{2}$ in the  plane
 ${\mathbb R}^{2}$, where $I$ is an open interval, is called {\it  convex} if, for all $s\in I$
 the trace $X(I)$ of $X$ lies entirely on one side of the closed
 half-plane determined by the tangent line at $X(s)$ (\cite{dC}).
 A regular plane curve $X:I\rightarrow {\mathbb R}^{2}$ is called {\it locally  convex} if, for each $s\in I$
 there exists an open  subinterval $I_0\subset I$ containing $s$ such that the curve $X|_{I_0}$ restricted to $I_0$
 is a convex curve.
\vskip 0.3cm
Henceforth,  we will  say  that a locally convex curve $X$ in the  plane
 ${\mathbb R}^{2}$ is  {\it strictly locally convex} if the curve   is smooth
 (that is, of class $C^{(3)}$) and is of positive  curvature $\kappa$
 with respect to the unit normal $N$ pointing to the convex side.
 Hence, in this case we have $\kappa(s)=\left< X''(s), N(X(s))\right> >0$,
  where $X(s)$ is an arc-length parametrization of $X$.
\vskip 0.3cm
When  $f:I\rightarrow {\mathbb R}$ is a smooth function defined on an open interval $I$,
we will also say that $f$ is
{\it strictly convex} if the graph of $f$ has  positive curvature $\kappa$
with respect to the upward unit normal $N$. This condition is equivalent to the positivity of $f''(x)$ on $I$.
\vskip 0.3cm

Suppose that $X$ is a strictly locally convex  curve  in the plane
 ${\mathbb R}^{2}$ with the unit normal $N$ pointing to the convex side.
 For a fixed point $P \in X$, and for a sufficiently small $h>0$, we consider the  line $\ell$ passing through
 $P+hN(P)$ which is parallel to
 the tangent line $t$ of $X$ at $P$ and the points  $A$ and $B$  where the line $\ell$ intersects the curve $X$.
We  denote by  $t_1$, $t_2$ the tangent
lines of $X$ at the points $A,B$ and by $Q,A_1,B_1$
the intersection points $t_1\cap t_2$,  $t_1\cap t$, $t_2\cap t$, respectively.

We let $L_P(h)$, $g_P(h)$, $j_P(h)$ and $k_P(h)$ denote the length of the chord $AB$ ,
 the distance from the centroid $G$  of the section of $X$ cut off by $\ell$
 to the line $\ell$,
 the distance from the centroid $J$   of the triangle  $\bigtriangleup QAB$
 to the line $\ell$ and  the distance from the centroid $K$   of the triangle  $\bigtriangleup QA_1B_1$
 to the line $\ell$, respectively.
\vskip 0.3cm
 Now, we consider  $S_P(h)$ and  $T_P(h)$
defined by the area of the region bounded by the curve $X$ and chord $AB$,
the area  $|\bigtriangleup PAB|$ of  triangle  $\bigtriangleup PAB$, respectively.
Then, obviously  we have
\begin{equation}\tag{1.1}
   \begin{aligned}
 T_P(h)=\frac{1}{2}hL_P(h)
    \end{aligned}
   \end{equation}
   and we get (\cite{KK4})
  \begin{equation}\tag{1.2}
   \begin{aligned}
 \frac{d}{dh}S_P(h)=L_P(h).
    \end{aligned}
   \end{equation}
\vskip 0.3cm
It is well known that  parabolas satisfy the following properties.

\vskip 0.3cm
\noindent {\bf Proposition 1.} Suppose that $X$ is an open part of a parabola.
 Then we have the following.

\noindent 1)  For arbitrary point $P\in X$ and sufficiently small $h>0$, $X$ satisfies
  \begin{equation}\tag{1.3}
   \begin{aligned}
  S_P(h)=\frac{4}{3}T_P(h).
    \end{aligned}
   \end{equation}

 \noindent 2) For arbitrary point $P\in X$ and sufficiently small $h>0$, $X$ satisfies
  \begin{equation}\tag{1.4}
   \begin{aligned}
  g_P(h)=\frac{2}{5}h.
    \end{aligned}
   \end{equation}

 \noindent 3)  For arbitrary point $P\in X$ and sufficiently small $h>0$, $X$ satisfies
   \begin{equation}\tag{1.5}
   \begin{aligned}
j_P(h)=\frac{2}{3}h.
    \end{aligned}
   \end{equation}

 \noindent 4)  For arbitrary point $P\in X$ and sufficiently small $h>0$, $X$ satisfies
   \begin{equation}\tag{1.6}
   \begin{aligned}
k_P(h)=\frac{4}{3}h.
    \end{aligned}
   \end{equation}

  \vskip 0.3cm
 \noindent {\bf Proof.}  For a proof of 1), see \cite{S}.
 If we denote by $V$ the point where the parallel line $m$ through
 the point $P$  to the axis of $X$ meets the chord $AB$, then $V$ is the mid point of
 $AB$ and the point $Q$ is on the line $m$ with $PV=PQ$.
 This completes the proof of 2), 3) and 4).   $\square$

 \vskip 0.3cm
 Very recently, the first author  of the present paper and Y. H. Kim
  showed that among strictly convex plane curves,
the above area property (1.3) of parabolic sections characterize
parabolas. More precisely, they proved as follows (\cite{KK4}).
\vskip 0.3cm

  \noindent {\bf Proposition 2.}
 Let $X$ be a strictly convex curve  in the plane
 ${\mathbb R}^{2}$. Then $X$ is a parabola if and only if  it satisfies
\vskip 0.3cm
 \noindent $(C):$
 For a point $P$ on $X$ and a chord $AB$ of $X$ parallel to the tangent of $X$ at $P$,
 the area of the region bounded by the curve and $AB$ is $4/3$ times
the area of triangle $\bigtriangleup ABP$, that is,
 $$
  S_P(h)=\frac{4}{3}T_P(h).
   $$

 \vskip 0.3cm

  Archimedes showed that parabolas satisfy (1.3) (\cite{S}).
 Actually, in \cite{KK4} the first author  of the present paper with Y. H. Kim
  established five  characterizations of parabolas, which are the converses of
well-known properties of parabolas  originally due to Archimedes
(\cite{S}).  For some properties and characterizations of parabolas with respect to the area of triangles associated with a curve,
see \cite{D,KKKP,KS, Kr}.
For the higher dimensional analogues of some results in \cite{KK4}, see \cite{KK2} and \cite{KK3}.

\vskip 0.3cm
 \vskip 0.3cm
In \cite{KKP}, using Proposition 1, D.-S. Kim et al.  proved the following characterization theorem for parabolas
with respect to the function $ g_P(h)$.
 \vskip 0.3cm
 \noindent {\bf Proposition 3.}
 Let $X$ be  a strictly locally convex plane curve  in the  plane
 ${\mathbb R}^{2}$. For a fixed point $P$ on $X$ and a sufficiently small $h>0$,
 we denote by $\ell$ the parallel line through $P+hN(P)$ to the tangent $t$ of the curve $X$ at $P$.
 We let $g_P(h)$ the distance from the center $G$ of gravity  of the section of $X$ cut off by $\ell$
 to the line $\ell$. Then $X$ is an open part of a parabola if and only if
 it satisfies
  for a fixed point $P$ on $X$ and a sufficiently small $h>0$
 $$
   g_P(h)=\frac{2}{5}h.
  $$
 \vskip 0.3cm
In \cite{KKP}, the distance from the center $G$ of gravity  of the section of $X$ cut off by $\ell$
 to the tangent $t$ of $X$ at $P$ was denoted by $d_P(h)$. Hence, we see that $d_P(h)+g_P(h)=h$.
\vskip 0.3cm
 In this article,  we
study whether the remaining properties of parabolas in Proposition 1  characterize parabolas.
 \vskip 0.3cm
First of all, in Section 2 we prove the following:
 \vskip 0.3cm

 \noindent {\bf Theorem 4.}
 Suppose that  $X$ denotes  a strictly locally convex plane curve  in the  plane
 ${\mathbb R}^{2}$. For a fixed point $P$ on $X$ and a sufficiently small $h>0$,
 we denote by $\ell$ the parallel line through $P+hN(P)$ to the tangent $t$ of the curve $X$ at $P$.
 Then we have
   \begin{equation}\tag{1.7}
   \begin{aligned}
  \lim _{h\rightarrow 0}\frac{j_P(h)}{h}=\frac{2}{3}.
       \end{aligned}
   \end{equation}
   and
    \begin{equation}\tag{1.8}
   \begin{aligned}
  \lim _{h\rightarrow 0}\frac{k_P(h)}{h}=\frac{4}{3}.
       \end{aligned}
   \end{equation}
 \vskip 0.3cm

 Finally, with the aid of the characterization theorem of parabolas (Theorem 3 in \cite{KK4}),
  in Section 4 we prove the following.
\vskip 0.3cm
 \noindent {\bf Theorem 5.}
 Suppose that  $X$ denotes a strictly locally convex $C^{(3)}$ curve  in the plane
 ${\mathbb R}^{2}$. Then the following are equivalent.

\noindent 1) For all $P\in X$ and sufficiently small $h>0$, $X$ satisfies
 \begin{equation}\tag{1.9}
   \begin{aligned}
j_P(h)=\lambda(P)h^{\mu(P)},
    \end{aligned}
   \end{equation}
where $\lambda(P)$ and $\mu(P)$ are some functions.

\noindent 2) For all $P\in X$ and sufficiently small $h>0$, $X$ satisfies
$$
j_P(h)=\frac{2}{3}h.
   $$
\noindent 3) $X$ is an open part of a parabola.
\vskip 0.3cm
For the function $k_P(h)$, the similar argument as in the proof of Theorem 5 yields the following.
\vskip 0.3cm
 \noindent {\bf Theorem 6.}
 Suppose that  $X$ denotes a strictly locally convex $C^{(3)}$ curve  in the plane
 ${\mathbb R}^{2}$. Then the following are equivalent.

\noindent 1) For all $P\in X$ and sufficiently small $h>0$, $X$ satisfies
 \begin{equation}\tag{1.10}
   \begin{aligned}
k_P(h)=\lambda(P)h^{\mu(P)},
    \end{aligned}
   \end{equation}
where $\lambda(P)$ and $\mu(P)$ are some functions.

\noindent 2) For all $P\in X$ and sufficiently small $h>0$, $X$ satisfies
$$
k_P(h)=\frac{4}{3}h.
   $$
\noindent 3) $X$ is an open part of a parabola.
 \vskip 0.3cm
 \noindent {\bf Remark.} If we consider the distance $\delta_P(h)$ from the centroid of the triangle $\Delta PAB$
 to the parallel line $\ell$ through the point $P+hN(P)$ to the tangent $t$ of a  strictly locally convex  curve $X$ at $P$,
 then $X$ always satisfies $\delta_P(h)=\frac{1}{3}h$ for all sufficiently small $h>0$.

 \vskip 0.3cm
For some characterizations of parabolas or  conic sections by  properties of tangent lines, see
\cite{KKa} and \cite{KKPj}. In \cite{KK1}, using curvature function $\kappa$ and support function $h$
of a plane curve,
the second  author of the present paper and Y. H. Kim gave a
characterization of ellipses and hyperbolas centered at the origin.

 Among  the  graphs of functions, \'A. B\'enyi et al. proved some characterizations of parabolas (\cite{BSV1, BSV2}). In \cite{R}, B.  Richmond and T. Richmond  established a dozen necessary and sufficient conditions
 for the graph of a function to
be a parabola  by using elementary techniques.
\vskip 0.3cm
  Throughout this article, all curves are of class $C^{(3)}$ and connected, unless otherwise mentioned.
  \vskip 0.50cm

\section{Proof of Theorem 4}

\vskip 0.5cm
In this section, we prove Theorem 4.
First of all,
we  need the following lemma (\cite{KK4}) which is useful in the proof of main theorems.
\vskip 0.3cm

 \noindent {\bf Lemma 7.} Suppose that   $X$  is a strictly locally convex  $C^{(3)}$ curve  in the plane
 ${\mathbb R}^{2}$ with  the unit normal $N$ pointing to the convex side.  Then
 we have
  \begin{equation}\tag{2.1}
   \begin{aligned}
   \lim_{h\rightarrow 0} \frac{1}{\sqrt{h}}L_P(h)= \frac{2\sqrt{2}}{\sqrt{\kappa(P)}},
    \end{aligned}
   \end{equation}
 where $\kappa(P)$ is the curvature of $X$ at $P$ with respect to  the unit normal $N$ pointing to the convex side.
 \vskip 0.3cm

Now, we prove Theorem 4 as follows.

 Let us denote by  $X$  a strictly locally convex  $C^{(3)}$ curve in the Euclidean plane ${\mathbb R^2}$.
We  fix an arbitrary  point $P$ on $X$.
Then, we may take a coordinate system $(x,y)$
 of  ${\mathbb R}^{2}$ such that  $P$ is  the origin $(0,0)$ and $x$-axis is the tangent line $t$ of $X$ at $P$.
 Furthermore, we may regard $X$ to be locally  the graph of a non-negative strictly convex  function $f: {\mathbb R}\rightarrow {\mathbb R}$
 with $f(0)=f'(0)=0$. Then $N$ is the upward unit normal.

 Since the curve $X$ is of class $C^{(3)}$, the Taylor's formula of $f(x)$ is given by
 \begin{equation}\tag{2.2}
 f(x)= ax^2 + f_3(x),
  \end{equation}
where  $2a=f''(0)$ and $f_3(x)$ is an $O(|x|^3)$  function.
Noting that the curvature $\k$ of $X$ at $P$ is given by $\kappa(P)=f''(0)>0$, we see that $a$ is positive.

For a sufficiently small $h>0$, the line $\ell$ through $P+hN(P)$ and orthogonal to $N(P)$ is
given by $y=h$. We denote by $A(s,f(s))$ and $B(t,f(t))$
 the points where the line $\ell:y=h$ meets the curve $X$ with $s<0<t$.
 Then we have $f(s)=f(t)=h$.  The tangent lines $t_1$ and $t_2$ to $X$ at $A$ and $B$
  intersect at
 the point $Q=(x_0(h),y_0(h))$ with
  \begin{equation}\tag{2.3}
x_0(h)=\frac{tf'(t)-sf'(s)}{ f'(t)-f'(s)},
  \end{equation}
  \begin{equation}\tag{2.4}
y_0(h)=h+\frac{(t-s)f'(t)f'(s)}{ f'(t)-f'(s)}<0
  \end{equation}
and they
 meet the $x$-axis (the tangent to $X$ at $P$) at
 $B_1(s-h/f'(s),0)$ and  $B_2(t-h/f'(t),0)$, respectively.

Noting $L_P(h)=t-s$, one gets
\begin{equation}\tag{2.5}
j_P(h)=h-\frac{1}{3}\{2h+y_0(h)\}=-\frac{1}{3}\frac{L_P(h)f'(t)f'(s)}{ f'(t)-f'(s)}
  \end{equation}
  and
  \begin{equation}\tag{2.6}
k_P(h)=h-\frac{1}{3}y_0(h)=\frac{2}{3}h-\frac{1}{3}\frac{L_P(h)f'(t)f'(s)}{ f'(t)-f'(s)}.
  \end{equation}
 Hence we  obtain
 \begin{equation}\tag{2.7}
\frac{j_P(h)}{h}=\frac{1}{3}\frac{L_P(h)}{\sqrt{h}}\frac{1}{\alpha_P(h)}
  \end{equation}
  and
  \begin{equation}\tag{2.8}
\frac{k_P(h)}{h}=\frac{2}{3}+\frac{j_P(h)}{h},
  \end{equation}
 where we use
 \begin{equation}\tag{2.9}
  \begin{aligned}
\alpha_P(h)=\frac{(f'(s)-f'(t))}{f'(s)f'(t)}\sqrt{h}.
   \end{aligned}
  \end{equation}

On the other hand, it follows from  Lemma 5 in \cite{KS} that
  \begin{equation}\tag{2.10}
  \begin{aligned}
\lim_{h\rightarrow0}\alpha_P(h)=\frac{\sqrt{2}}{\sqrt{\kappa(P)}}.
   \end{aligned}
  \end{equation}
Together with (2.7) and Lemma 7, this shows that
\begin{equation}\tag{2.11}
  \begin{aligned}
\lim_{h\rightarrow0}\frac{j_P(h)}{h}=\frac{2}{3},
   \end{aligned}
  \end{equation}
and hence from (2.8) we also get
\begin{equation}\tag{2.12}
  \begin{aligned}
\lim_{h\rightarrow0}\frac{k_P(h)}{h}=\frac{4}{3}.
   \end{aligned}
  \end{equation}
This completes the proof of Theorem 4.
  \vskip 0.3cm

  \vskip 0.5cm
\section{Proofs of Theorems 5 and 6}

\vskip 0.5cm

In this section, in order to prove Theorem 5 we use the main result of \cite{KK4} (Theorem 3 in \cite{KK4})
and Theorem 4 stated in  Section 1.

First, we prove

 \vskip 0.3cm
\noindent {\bf Lemma 8.} Suppose that   $X$  is a strictly locally convex  $C^{(3)}$ curve  in the plane
 ${\mathbb R}^{2}$ with  the unit normal $N$ pointing to the convex side.  Then
 we have
 \begin{equation}\tag{3.1}
  \begin{aligned}
\sqrt{h}\frac{d}{dh}L_P(h)=\alpha_P(h).
   \end{aligned}
  \end{equation}
where $\alpha_P(h)$ is defined in (2.9).

   \vskip 0.3cm
 \noindent {\bf Proof.}
 Just as in the proof of Theorem 4 in Section 2, for an arbitrary  point $P$ on $X$ we take a coordinate system $(x,y)$
 of  ${\mathbb R}^{2}$ so that (2.2) holds with $f(0)=f'(0)=0$ and $2a=f''(0)>0$.
 Then, for sufficiently small $h>0$, we put $f(t)=h$ with $t>0$ and we denote by $A(s(t),h)$ and $B(t,h)$
 the points where the line $\ell:y=h$ meets the curve $X$ with $s=s(t)<0<t$.
 Then we have
  \begin{equation}\tag{3.2}
  \begin{aligned}
f(s(t))=f(t)=h
   \end{aligned}
  \end{equation}
  and
   \begin{equation}\tag{3.3}
  \begin{aligned}
L_P(h)=t-s(t).
   \end{aligned}
  \end{equation}

Noting $h=f(t)$, one obtains from (3.3) that
 \begin{equation}\tag{3.4}
  \begin{aligned}
\frac{d}{dh}L_P(h)=(1-s'(t))\frac{dt}{dh}=\frac{1-s'(t)}{f'(t)}.
   \end{aligned}
  \end{equation}
Therefore, it follows from (3.2) that
\begin{equation}\tag{3.5}
  \begin{aligned}
\frac{d}{dh}L_P(h)=\frac{1}{f'(t)}-\frac{1}{f'(s(t))}=\frac{f'(s)-f'(t)}{f'(t)f'(s)}.
   \end{aligned}
  \end{equation}
Together with (2.9), this  completes the proof of Lemma 8. $\square$
\vskip 0.3cm
 It is obvious that any open part of parabolas satisfy 1) and 2) in Theorem 5.

Conversely, suppose that $X$ is a strictly locally convex $C^{(3)}$ curve  in the plane
 ${\mathbb R}^{2}$ which  satisfies for all $P\in X$ and sufficiently small $h>0$
 $$
 j_P(h)=\lambda(P)h^{\mu(P)},
       $$
   where $\lambda(P)$ and $\mu(P)$ are some functions.
 Using Theorem 4, by letting $h\rightarrow 0$ we see that
 \begin{equation}\tag{3.6}
   \begin{aligned}
 \lim_{h\rightarrow 0}h^{\mu(P)-1}=\frac{2}{3}\frac{1}{\lambda(P)},
         \end{aligned}
   \end{equation}
which shows that
   $\mu(P)=1$ and   $\lambda(P)=\frac{2}{3}$.
Therefore, the curve $X$ satisfies for all $P\in X$ and sufficiently small $h>0$
$$
 j_P(h)=\frac{2}{3}h.
       $$

Now, using Lemma 8 we get the following.
\vskip 0.3cm
\noindent {\bf Lemma 9.} Suppose that  $X$ denotes a strictly locally convex $C^{(3)}$ curve  in the plane
 ${\mathbb R}^{2}$ which satisfies  (1.5) for all $P\in X$ and sufficiently small $h>0$.
Then for all $P\in X$ and sufficiently small $h>0$ we have

 \begin{equation}\tag{3.7}
 L_P(h)=\frac{2\sqrt{2}}{\sqrt{\kappa(P)}}\sqrt{h}.
  \end{equation}
\vskip 0.3cm

  \noindent {\bf Proof.} It follows from (1.5) and (2.7)  that
  \begin{equation}\tag{3.8}
 L_P(h)=2\sqrt{h}\alpha_P(h).
  \end{equation}
  Together with Lemma 8, this yields
 \begin{equation}\tag{3.9}
 2h\frac{d}{dh}L_P(h)=L_P(h).
  \end{equation}

By integrating (3.9) with respect to $h$, we get for some constant $C=C(P)$
  \begin{equation}\tag{3.10}
L_P(h)=C\sqrt{h}.
  \end{equation}
Thus, Lemma 7 completes the proof of Lemma 9. $\square$
 \vskip 0.3cm

\vskip 0.3cm

Finally, we prove Theorem 5 as follows.
\vskip 0.3cm
It follows from (1.2) and  $S_P(0)=0$ that
by integrating (3.7) we get
 \begin{equation}\tag{3.11}
 S_P(h)=\frac{4\sqrt{2}}{3\sqrt{\kappa(P)}}h\sqrt{h}.
  \end{equation}
 Hence, together with  (1.1) and (3.7), (3.11) yields that for all $P\in X$ and sufficiently small $h>0$
\begin{equation}\tag{3.12}
 S_P(h)=\frac{4}{3}T_P(h).
  \end{equation}
Thus, it follows from Proposition 2 that  $X$ is an open part of a parabola.
This  completes the proof of Theorem 5.
\vskip 0.3cm
In order to prove Theorem 6, with the help of (2.8) we may use the similar argument as in the proof of Theorem 5.
\vskip 0.50cm


\vskip 1.0 cm

Department of Mathematics, \par Chonnam National University,\par
Kwangju 500-757, Korea

{\tt E-mail: dongseo@chonnam.ac.kr}
\vskip 0.3 cm

Department of Mathematics, \par Chonnam National University,\par
Kwangju 500-757, Korea

{\tt E-mail: dosokim@chonnam.ac.kr}
\vskip 0.3 cm

 \vskip 0.50cm

\vskip 0.3 cm
\end{document}